\documentclass[11pt,a4paper]{article}

\usepackage[utf8]{inputenc}
\usepackage[T1]{fontenc}
\usepackage{lmodern}
\usepackage{amsmath,amssymb,amsthm,mathtools}
\usepackage{microtype}
\usepackage[margin=1in]{geometry}
\usepackage{authblk}
\usepackage{hyperref}
\hypersetup{
    colorlinks=true,
    linkcolor=blue!50!black,
    citecolor=blue!50!black,
    urlcolor=blue!50!black
}
\usepackage{tikz}
\usetikzlibrary{decorations.markings,decorations.pathmorphing,arrows.meta}

\definecolor{blueamu}{RGB}{30,80,160}
\definecolor{darkblueamu}{RGB}{15,40,95}

\newtheorem{theorem}{Theorem}

\theoremstyle{remark}

\newcommand{\E}{\mathbb{E}}
\newcommand{\R}{\mathbb{R}}
\newcommand{\Li}{\operatorname{Li}}
\newcommand{\Res}{\operatorname{Res}}
\DeclareMathOperator{\Real}{Re}
\DeclareMathOperator{\Imag}{Im}
\newcommand{\dd}{\mathrm{d}} 

\title{Closed-form linear moments of the two-dimensional\\
angular central Gaussian distribution}

\author{ Siméon Vareilles}
\affil{Aix Marseille Univ, Universit\'e de Toulon, CNRS, CPT, Marseille, France\\\texttt{vareilles@cpt.univ-mrs.fr}}

\date{\today}

\begin{document}
\maketitle

\begin{abstract}
The polar-angle marginal of a centred bivariate Gaussian distribution, obtained after integrating out the radial coordinate, gives the two-dimensional angular central Gaussian (ACG) distribution of Tyler.  While its trigonometric and vector-valued moments have been studied in detail, to our knowledge there are no explicit closed-form expressions for the \emph{linear} moments $\E[\theta]$ and $\E[\theta^{2}]$ on the natural domain $\theta\in\left]-\pi/2,\pi/2\right[$.  Here \emph{linear} refers to the ordinary moments $\int\theta^{k}f(\theta)\,d\theta$ of the angle regarded as a real-valued variable, in contrast to the circular (trigonometric) moments $\E[e^{ik\theta}]$ customary in directional statistics.  We provide such expressions: the mean is a simple arctangent of the parameters, while the second moment is given by the real part of a dilogarithm.  The derivation, based on a contour integration around the branch cut of $\arctan z$, is elementary.  These quantities naturally arise in physics, where $\theta$ is interpreted as a real-valued phase rather than a circular variable.
\end{abstract}

\section{Introduction}
\label{sec:intro}

Centred Gaussian probability distributions on $\R^{2}$ are ubiquitous in physics. In classical statistical mechanics they describe equilibrium fluctuations of conjugate variables; in signal and image processing they arise as noise models on the Stokes plane; in continuous-variable quantum information they are the Wigner functions of zero-displacement Gaussian states, which form the primary class of bosonic states used in quantum optics and quantum computing or in cosmology to study the property of the Gaussian state of perturbations
\cite{weedbrook2012gaussian,adesso2014continuous,serafini2023quantum}. A single-mode Gaussian state is then fully
specified by its $2\times 2$ covariance matrix $\Sigma$ in the
canonical pair $(q,p)\in\R^{2}$ of \emph{position} and
\emph{momentum}-like variables.

Introducing polar coordinates,
\begin{equation}
    q = r\cos\theta,
    \qquad
    p = r\sin\theta,
    \qquad r\ge 0,\quad
    \theta\in\,\left]-\tfrac{\pi}{2},\tfrac{\pi}{2}\right[,
\end{equation}
where $r\ge 0$ and the half-open range $\theta\in\left]-\pi/2,\pi/2\right[$ gives a
one-to-one parametrisation of the oriented lines through the origin, reflecting the
identification $(r,\theta)\sim(-r,\theta+\pi)$. Integrating out $r$ yields the marginal
probability density function (pdf) of the \emph{polar angle},
\begin{equation}\label{eq:f}
    f(\theta)
    =\frac{1}{2\pi}\,
    \frac{1}{A\cos^{2}\theta-2C\sin\theta\cos\theta+B\sin^{2}\theta},
    \qquad \theta\in\bigl]-\tfrac{\pi}{2},\tfrac{\pi}{2}\bigr[,
\end{equation}
where the boundary $\theta=\pi/2$ carries zero probability mass.\footnote{For a general centred bivariate Gaussian the normalised angular density on the
half-circle $\left]-\pi/2,\pi/2\right[$ is
$f(\theta)=\dfrac{\sqrt{AB-C^{2}}}{\pi}\,
\big(A\cos^{2}\theta-2C\sin\theta\cos\theta+B\sin^{2}\theta\big)^{-1}$.
The simple prefactor $1/(2\pi)$ in \eqref{eq:f} corresponds to the pure single-mode case,
for which the Heisenberg uncertainty relation is saturated and
$\det\Sigma=AB-C^{2}=1/4$ (using $[\hat q,\hat p]=i$, so that the vacuum has
$\Sigma=\tfrac12 I$).}

The real parameters $A,B,C$ are the second moments of the centred Gaussian state, i.e.\ the
entries of the covariance matrix $\Sigma$,\footnote{Explicitly $A=\langle p^{2}\rangle$,
$B=\langle q^{2}\rangle$, $C=\langle qp\rangle$. Note that $A=\langle p^{2}\rangle$ is the
coefficient of $\cos^{2}\theta$ even though $q=r\cos\theta$: the angular weight is
$\propto(\hat n^{\top}\Sigma^{-1}\hat n)^{-1}$ with $\hat n=(\cos\theta,\sin\theta)$, and the
$(1,1)$ entry of $\Sigma^{-1}$ is $\langle p^{2}\rangle/\det\Sigma$. Equivalently, a larger
$\langle p^{2}\rangle$ concentrates the orientation toward the $p$-axis ($\theta=\pi/2$).}
and satisfy
\begin{equation}\label{eq:constraint}
    A>0,\qquad B>0,\qquad AB-C^{2}=\det\Sigma=\tfrac14\,.
\end{equation}
The pdf~\eqref{eq:f} is $\pi$-periodic and is therefore most naturally defined on the
projective line, here identified with $\left]-\pi/2,\pi/2\right[$.

Distribution~\eqref{eq:f} is far from new.  It is precisely the two-dimensional case of the \emph{angular central Gaussian (ACG) distribution}, originally introduced by Tyler \cite{tyler1987statistical} as an alternative to the Bingham distribution for antipodally symmetric directional data, and further discussed in Mardia and Jupp textbook~\cite{mardia2009directional}. Under the angle-doubling map $\varphi=2\theta$, this density function transforms into the \emph{wrapped Cauchy} distribution on the full circle \cite{kent2022directional,kato2013extended}.  Moments of \emph{vector-valued} quantities such as $\E[\mathbf{p}\otimes\mathbf{p}]$ (the so-called ``orientation tensor'') are commonly used in directional statistics \cite{tyler1987statistical,mardia2009directional}. Additionally, the \emph{trigonometric} moments $\E[\mathrm{e}^{ik\theta}]$ have an analytical expression inherited from the wrapped-Cauchy structure \cite{mardia2009directional,kato2013extended}.

In contrast, the natural \emph{linear} moments
\begin{equation}\label{eq:moments_def}
    \E[\theta^{n}]
    =\int_{-\pi/2}^{\pi/2}\theta^{n}\,f(\theta)\,\dd\theta,
\end{equation}
defined on the fundamental domain $\left]-\pi/2,\pi/2\right[$, have received much less attention.  In the context of circular-statistics, this is expected because a linear moment depends on the choice of branch and is not invariant under rotations of the underlying circle.  In physics, however, $\theta$ is often a genuinely real variable with direct physical significance, such as a quadrature
angle, a homodyne phase, or the orientation of a polarisation ellipse. In these cases, the variance of $\theta$ serves as the relevant figure of merit.  A closely related non-centred isotropic case has been considered in the context of the angle of linear polarisation \cite{kupinski2014relating}, although in that work the underlying Gaussian has nonzero mean and therefore leads to a different distribution.

In this note we provide closed-form expressions for the first two
linear moments of~\eqref{eq:f}.  The result is:

\begin{theorem}\label{thm:main}
    Let $f(\theta)$ be the pdf defined by~\eqref{eq:f}--%
    \eqref{eq:constraint}.  Then
    \begin{align}
        \E[\theta]   &= \arctan\!\frac{2C}{2B+1},
        \label{eq:mean}\\[4pt]
        \E[\theta^{2}] &= \frac{\pi^{2}}{12}
        +\Real\!\left[\,
        \Li_{2}\!\left(1-\frac{4B}{1+2B-2\,i\,C}\right)\,\right],
        \label{eq:var}
    \end{align}
    where $\Li_{2}(z)=\sum_{k\ge 1}z^{k}/k^{2}$ is the dilogarithm
    function.  The variance is then $\operatorname{Var}(\theta)=\E[\theta^{2}]-\E[\theta]^{2}$.
\end{theorem}

Formulas~\eqref{eq:mean}--\eqref{eq:var} constitute the main contribution of this paper.  The remainder is structured as follows. Section~\ref{sec:setup} introduces the contour-integral representation of the moments, Section~\ref{sec:derivation} provides the derivation of Theorem~\ref{thm:main}, and Section~\ref{sec:checks} discusses limiting cases and conducts brief consistency checks.

\section{Contour-integral representation}\label{sec:setup}

The constraints~\eqref{eq:constraint} ensure that the denominator of
$f(\theta)$ has no real zeros and that $f$ integrates to one.
Indeed, the change of variables
\begin{equation}
    z=\tan\theta,\qquad
    z\in\left]-\infty,\infty\right[,\qquad
    \dd\theta=\frac{\dd z}{1+z^{2}},
\end{equation}
maps the denominator of~\eqref{eq:f} to the quadratic
$A-2Cz+Bz^{2}$ (multiplied by $\cos^{2}\theta$, which cancels
$1/(1+z^{2})$).  Hence the moments take the compact form
\begin{equation}\label{eq:contour_form}
    \E[\theta^{n}]
    =\frac{1}{2\pi B}\int_{-\infty}^{\infty}
    \frac{\arctan^{n}z}{(z-z_{+})(z-z_{-})}\,\dd z,
\end{equation}
where the poles are complex conjugates of one another,
\begin{equation}
    z_{\pm}=\frac{1}{B}\!\left(C\pm\frac{i}{2}\right),
    \qquad
    z_{+}-z_{-}=\frac{i}{B}.
\end{equation}

Throughout, we take the principal branches of both $\ln$ and $\Li_{2}$.

\begin{figure}[htbp]
    \centering
    \begin{tikzpicture}[
        scale=1.3,
        >=Stealth,
        line cap=round,
        every node/.style={font=\small}
    ]
    \draw[->,black] (-3.5,0) -- (3.7,0) node[below right,black] {$\Real z$};
    \draw[->,black] (0,-1)  -- (0,3.8) node[left,black]      {$\Imag z$};
    \draw[purple,thick,decorate,
          decoration={zigzag,segment length=4pt,amplitude=2pt}]
        (0,1) -- (0,3.5);
    \node[purple,anchor=west,font=\small\itshape] at (0.10,3.2) {branch cut};
    \filldraw[fill=white,draw=black,line width=0.5pt] (0,1) circle (1.6pt);
    \node[anchor=west] at (0.20,1.06) {$+i$};
    \filldraw[darkblueamu] (1.6, 0.35) circle (1.8pt)
        node[anchor=west,blue!40!black,xshift=2pt] {$z_+$};
    \filldraw[darkblueamu] (1.6,-0.35) circle (1.8pt)
        node[anchor=west,blue!40!black,xshift=2pt] {$z_-$};
    \node[below=2pt] at (-3,0) {$-R$};
    \node[below=2pt] at ( 3,0) {$R$};
    \draw[blueamu,very thick,
          decoration={markings,
            mark=at position 0.15  with {\arrow{>}},
            mark=at position 0.42  with {\arrow{>}},
            mark=at position 0.59  with {\arrow{>}},
            mark=at position 0.655 with {\arrow{>}},
            mark=at position 0.72  with {\arrow{>}},
            mark=at position 0.88  with {\arrow{>}}
          },
          postaction={decorate}]
        (-3,0) -- (3,0)
        arc[start angle=0,end angle=87,radius=3]
        -- (0.157,1)
        arc[start angle=0,end angle=-180,radius=0.157]
        -- (-0.157,2.996)
        arc[start angle=93,end angle=180,radius=3];
    \draw[gray!70,thin,->,shorten >=1.5pt] (-0.55,0.55) -- (-0.10,0.92);
    \node[gray] at (-0.65,0.48) {$\epsilon$};
    \draw[gray!80,dashed] (0,0) -- ({3*cos(45)},{3*sin(45)});
    \node[gray] at ({1.6*cos(45)+0.18},{1.6*sin(45)-0.05}) {$R$};
    \node[blueamu,font=\large] at (-2.3,2.2) {$\gamma$};
    \end{tikzpicture}
    \caption{Integration contour $\gamma$ used to
    evaluate~\eqref{eq:contour_form}.  The contour closes in the
    upper half-plane, encircles the pole $z_{+}$ and wraps the
    branch cut of $\arctan z$ along $\{it:t\ge 1\}$.  In the limits
    $R\to\infty$ and $\epsilon\to 0$, the contributions of the
    large arc and of the small indentation around $z=i$ vanish.}
    \label{fig:contour}
\end{figure}
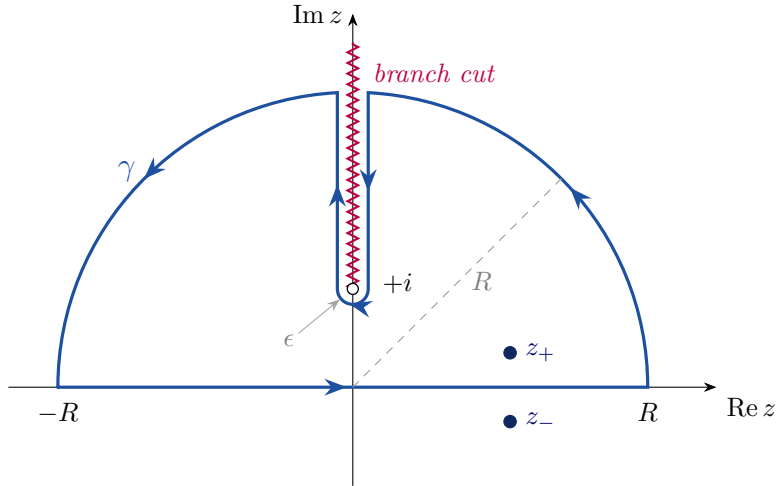

For the principal branch
\begin{equation}\label{eq:arctan_log}
    \arctan z
    =\frac{1}{2i}\ln\!\frac{1+iz}{1-iz},
\end{equation}
the function $\arctan z$ has two branch cuts along the imaginary
axis, $z=\pm it$ with $t\ge 1$.  We close the integration
contour~$\gamma$ by a large semicircle in the upper half-plane and
detour around the cut on the positive imaginary axis, see
Fig.~\ref{fig:contour}.  The only pole enclosed is $z_{+}$, and
both the large arc and the small indentation around the branch
point $z=i$ contribute negligibly in the limit $R\to\infty$,
$\epsilon\to 0$ (the latter because $\arctan z$ diverges only
logarithmically near $z=i$, whereas $\dd z=O(\epsilon)$).

\section{Derivation}\label{sec:derivation}

\subsection{Zeroth moment (normalisation)}

For $n=0$ the integrand has no branch cut, and the contour reduces
to the real axis closed by the large semicircle.  The integrand
decays as $|z|^{-2}$, so the arc contribution vanishes.  Cauchy's
theorem gives
\begin{equation}
    \E[1]
    =\frac{1}{2\pi B}\cdot 2\pi i\,
    \Res_{z=z_{+}}\!\frac{1}{(z-z_{+})(z-z_{-})}
    =\frac{1}{2\pi B}\cdot\frac{2\pi i}{z_{+}-z_{-}}
    =1,
\end{equation}
confirming that~\eqref{eq:f} is correctly normalised under
constraint~\eqref{eq:constraint}.

\subsection{First moment}

For $n=1$ the residue contribution is
\begin{equation}
    2\pi i\,\Res_{z=z_{+}}\!\frac{\arctan z}{(z-z_{+})(z-z_{-})}
    =\frac{2\pi i\,\arctan z_{+}}{z_{+}-z_{-}}
    =2\pi B\,\arctan z_{+}.
\end{equation}
Across the upper cut the principal branch jumps by
$\arctan(it+0)-\arctan(it-0)=\pi$ ($t>1$).  Since the contour
descends the right edge of the cut and ascends the left edge (see
Fig.~\ref{fig:contour}), the discontinuity enters with a relative
minus sign, and the branch-cut contribution reads
\begin{align}
    \Delta_{\rm cut}
    &=-\,i\,\pi\int_{1}^{\infty}\frac{\dd t}{(it-z_{+})(it-z_{-})}
    \nonumber\\
    &=-\,i\pi B\,
    \ln\!\frac{i-z_{+}}{i-z_{-}},
\end{align}
where the last equality follows from a partial-fraction
decomposition and direct integration of $\int \dd t/(it-z_{\pm})$.

Combining the residue and cut contributions through
$\int_{-\infty}^{\infty}=2\pi i\,\Res_{z_{+}}-\Delta_{\rm cut}$ and
dividing by $2\pi B$,
\begin{equation}
    \E[\theta]
    =\arctan z_{+}+\frac{i}{2}\ln\!\frac{i-z_{+}}{i-z_{-}}.
\end{equation}
Substituting the explicit form of $z_{\pm}$ and simplifying gives
\begin{equation}
    \E[\theta]
    =-\frac{i}{2}\,\ln\!\frac{1+2B+2iC}{1+2B-2iC}
    =\arctan\!\frac{2C}{2B+1},
\end{equation}
which is the formula stated in Theorem~\ref{thm:main}.

\subsection{Second moment}

For $n=2$ we proceed in exactly the same way.  The residue
contribution is $2\pi B\,\arctan^{2}z_{+}$, and, with the same
orientation of the cut as above, the cut contribution is
\begin{equation}
    \Delta_{\rm cut}
    =i\int_{1}^{\infty}
    \frac{[\arctan(it-0)]^{2}-[\arctan(it+0)]^{2}}
    {(it-z_{+})(it-z_{-})}\,\dd t.
\end{equation}
Using the principal-branch representation
\begin{equation}
    \arctan(it\pm 0)=\pm\frac{\pi}{2}+\frac{1}{2i}\ln\!\frac{t-1}{t+1},
    \qquad t>1,
\end{equation}
the squared jump simplifies 
\begin{equation}
    [\arctan(it-0)]^{2}-[\arctan(it+0)]^{2}
    =i\pi\,\ln\!\frac{t-1}{t+1}.
\end{equation}
The cut contribution then becomes
\begin{equation}
    \Delta_{\rm cut}
    =-\pi\int_{1}^{\infty}
    \frac{\ln\!\dfrac{t-1}{t+1}}
    {(it-z_{+})(it-z_{-})}\,\dd t.
\end{equation}

We evaluate this integral by the substitution
\begin{equation}
    u=\frac{t-1}{t+1},
    \qquad t=\frac{1+u}{1-u},
    \qquad \dd t=\frac{2\,du}{(1-u)^{2}}.
\end{equation}
Each factor in the denominator transforms as
$it-z_{\pm}=\bigl[(i-z_{\pm})+u(i+z_{\pm})\bigr]/(1-u)$, so the
two factors of $(1-u)^{-2}$ cancel and
\begin{equation}
    \Delta_{\rm cut}
    =-2\pi\int_{0}^{1}
    \frac{\ln u}
    {\bigl[(i-z_{+})+u(i+z_{+})\bigr]
    \bigl[(i-z_{-})+u(i+z_{-})\bigr]}\,du.
\end{equation}
The roots in $u$ are
\begin{equation}
    \alpha_{\pm}=\frac{z_{\pm}-i}{z_{\pm}+i}.
\end{equation}
A partial-fraction split and use of the identity
$(i+z_{+})(i+z_{-})(\alpha_{+}-\alpha_{-})=2i(z_{+}-z_{-})=-2/B$
gives
\begin{equation}
    \Delta_{\rm cut}
    =\pi B\int_{0}^{1}
    \!\left(\frac{\ln u}{u-\alpha_{+}}
    -\frac{\ln u}{u-\alpha_{-}}\right)\dd u.
\end{equation}
Each integral is elementary (for $\alpha\notin[0,1]$),
\begin{equation}
    \int_{0}^{1}\frac{\ln u}{u-\alpha}\,\dd u
    =\Li_{2}\!\left(\frac{1}{\alpha}\right),
\end{equation}
(the logarithmic boundary term $\ln u\,\ln(1-u/\alpha)$ vanishes at
both endpoints, and $\Li_{2}(0)=0$).  Consequently
\begin{equation}
    \Delta_{\rm cut}
    =\pi B\!\left[\Li_{2}\!\left(\frac{1}{\alpha_{+}}\right)
    -\Li_{2}\!\left(\frac{1}{\alpha_{-}}\right)\right].
\end{equation}

Collecting all contributions through
$\E[\theta^{2}]=\arctan^{2}z_{+}-\Delta_{\rm cut}/(2\pi B)$,
\begin{equation}\label{eq:second_intermediate}
    \E[\theta^{2}]
    =\arctan^{2}z_{+}
    +\tfrac12
    \!\left[\Li_{2}\!\left(\frac{z_{-}+i}{z_{-}-i}\right)
    -\Li_{2}\!\left(\frac{z_{+}+i}{z_{+}-i}\right)\right].
\end{equation}

To bring~\eqref{eq:second_intermediate} into the form stated in
Theorem~\ref{thm:main}, introduce
\begin{equation}\label{eq:wdef}
    w=1-\frac{4B}{1+2B-2iC}=\frac{1-2B-2iC}{1+2B-2iC}.
\end{equation}
A direct calculation shows that
\begin{equation}
    \frac{z_{-}+i}{z_{-}-i}=w^{*},
    \qquad
    \frac{z_{+}+i}{z_{+}-i}=\frac{1}{w},
\end{equation}
so the two dilogarithms have arguments that are conjugate-and-inverse
to each other.  Using the conjugation identity
$\Li_{2}(w^{*})=\Li_{2}(w)^{*}$ and the inversion formula
\begin{equation}
    \Li_{2}\!\left(\frac{1}{w}\right)
    =-\frac{\pi^{2}}{6}-\frac12\ln^{2}(-w)-\Li_{2}(w),
\end{equation}
the bracket in~\eqref{eq:second_intermediate} becomes
\begin{equation}
    \tfrac12\!\left[\Li_{2}(w^{*})-\Li_{2}(1/w)\right]
    =\Real[\Li_{2}(w)]+\frac{\pi^{2}}{12}+\tfrac14\ln^{2}(-w).
\end{equation}
Finally, from \eqref{eq:arctan_log},
\begin{equation}
    \frac{1+iz_{+}}{1-iz_{+}}
    =\frac{-1+2B+2iC}{1+2B-2iC}=-w,
\end{equation}
so $\arctan^{2}z_{+}=-\tfrac14\ln^{2}(-w)$ and the
$\ln^{2}(-w)$ terms cancel exactly.  We arrive at
\begin{equation}
    \E[\theta^{2}]
    =\frac{\pi^{2}}{12}
    +\Real\!\left[\Li_{2}\!\left(
    1-\frac{4B}{1+2B-2iC}\right)\right],
\end{equation}
which completes the proof of Theorem~\ref{thm:main}.\qed

\section{Limiting cases and consistency checks}\label{sec:checks}

\paragraph{Isotropic case ($A=B=1/2$, $C=0$).}
The pdf reduces to the uniform distribution on $\left]-\pi/2,\pi/2\right[$.
Formula~\eqref{eq:mean} gives $\E[\theta]=\arctan 0=0$, and
$w=0$ in~\eqref{eq:wdef} so $\Li_{2}(0)=0$ and
$\E[\theta^{2}]=\pi^{2}/12$, in agreement with the elementary
variance of the uniform law on $\left]-\pi/2,\pi/2\right[$.

\paragraph{Concentrated limit along the $q$-axis ($B\to\infty$).}
The pdf becomes a Dirac peak at $\theta=0$.  In~\eqref{eq:wdef},
$w\to -1$, $\Li_{2}(-1)=-\pi^{2}/12$, and so
$\E[\theta^{2}]\to 0$, as expected.

\paragraph{Concentrated limit along the $p$-axis ($A\to\infty$).}
The constraint~\eqref{eq:constraint} forces $B\to 0$ with $C$
bounded, and the pdf becomes a sum of Dirac peaks at
$\theta=\pm\pi/2$.  Here $w\to 1$, $\Li_{2}(1)=\pi^{2}/6$, and so
$\E[\theta^{2}]\to\pi^{2}/12+\pi^{2}/6=\pi^{2}/4$, consistent with
$\theta=\pm\pi/2$ almost surely.

\paragraph{Reflection symmetry.}
Replacing $C\mapsto-C$ leaves the pdf invariant under
$\theta\mapsto-\theta$.  Formula~\eqref{eq:mean} flips sign as
required, and~\eqref{eq:var} is unchanged because $w$ becomes
$w^{*}$ and $\Real\Li_{2}(w^{*})=\Real\Li_{2}(w)$.

\paragraph{Connection with the wrapped Cauchy.}
The substitution $\varphi=2\theta$ maps $f(\theta)$ on
$\left]-\pi/2,\pi/2\right[$ to a wrapped Cauchy density on $\left]-\pi,\pi\right[$ with
concentration parameter $\rho$ and mean direction $\mu$ satisfying
$\rho\,\mathrm{e}^{i\mu}=-w=\dfrac{2B-1+2iC}{1+2B-2iC}$.  In particular,
$|w|<1$ is the wrapped-Cauchy mean-resultant length, a fact that
underlies the convergence of $\Li_{2}(w)$ in~\eqref{eq:var}.

\section*{Conclusion}

We have derived simple closed-form expressions for the linear mean and second moment of the two-dimensional angular central Gaussian distribution on the domain $\left]-\pi/2,\pi/2\right[$.  The mean is given by an arctangent, while the second moment involves the dilogarithm evaluated at the complex parameter
\begin{equation}\label{eq:wdef-concl}
    w=1-\frac{4B}{1+2B-2iC}\,,
\end{equation}
which is a fractional linear (Möbius) transformation of the covariance parameters.\footnote{In the present notation, this is the Möbius image of the pair \((B,C)\) that naturally appears after the contour calculation.}  These results complement the well-known formulas for the trigonometric and vector-valued moments of the 2D ACG / wrapped Cauchy family, and may be useful in physical applications where the polar angle of a Gaussian state is considered as a real-valued rather than a circular variable. To our knowledge, these closed-form linear moments do not appear explicitly in the directional-statistics literature.




\bibliographystyle{JHEP}
\bibliography{biblio-2}

\end{document}